\documentclass{ifacconf}
\usepackage{graphicx}      
\usepackage[square]{natbib}
\usepackage{amssymb}
\usepackage{amstext}
\usepackage{amsmath}
\usepackage{float}
\usepackage{color}

\begin{document}
\begin{frontmatter}

\title{Optical Costas loop: pull-in range estimation and hidden oscillations}

\author[spb,fin]{Kuznetsov N.~V.}
\author[spb,ipmash]{Leonov G.~A.}
\author[spb]{Seledzhi S.~M.}
\author[spb]{Yuldashev M.~V.}
\author[spb]{Yuldashev R.~V.}

\address[spb]{Faculty of Mathematics and Mechanics,
Saint-Petersburg State University, Russia}
\address[fin]{Dept. of Mathematical Information Technology,
University of Jyv\"{a}skyl\"{a}, Finland}
\address[ipmash]{Institute for Problems in Mechanical Engineering of
the Russian Academy of Sciences, Russia}

\begin{abstract}
In this work we consider a mathematical model of the optical Costas loop.
The pull-in range of the model is estimated by analytical and numerical methods.
Difficulties of numerical analysis, related to the existence of so-called hidden oscillations
in the phase space, are discussed.
\end{abstract}


\end{frontmatter}

\section{Introduction}
The Costas loop is a special mdification of the phase-locked loop,
which is widely used in telecommunication
for the data demodulation and carrier recovery.
The optical Costas loop is used
in intersatellite communication, see e.g. European Data Relay System (EDRS)
\citep{rosenkranz2016receiver,schaefer2015comparison,heine2014european}.
In the optical intersatellite communication
the lasers are used instead of radio signals
for the data transmission between satellites.
Both lasers, used in the transmitter and the receiver,
have a frequency mismatch due to the Doppler shift,
natural frequency drift, and phase noise.
The optical Costas loop adjusts the frequency and phase of the receiver
to the incoming signal,
which allows one to demodulate received data.
Operating principles of this circuit are similar to those
of conventional Costas loop circuit (homodyne detection).

In Section~2 a mathematical model of the optical Costas loop is derived.
In Section~3 the pull-in range of the optical Costas loop
with active proportional-integral (PI)
and lead-lag filters are discussed.

\section{Mathematical model of Optical Costas loop}
Consider a nonlinear mathematical model of the optical Costas loop model
in the signal space  (see, e.g. \citep{SchaeferR-2015})
\begin{figure}[H]
  \centering
  \includegraphics[width=\linewidth]{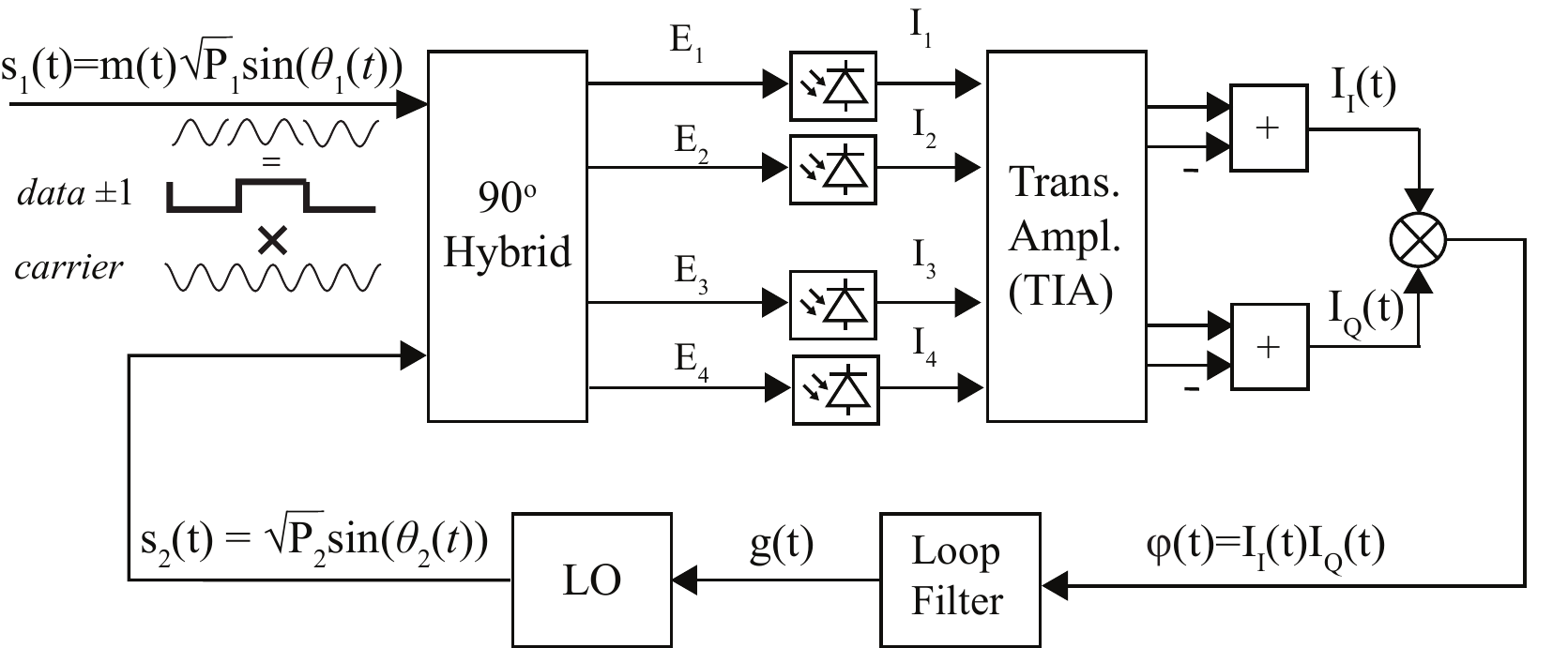}
  \caption{Optical BPSK Costas loop in the signal space}
\label{costas_after_sync}
\end{figure}
The input signal is a BPSK signal,
which is the product of the transferred data $m(t) = \pm 1$
and the harmonic high frequency carrier $\sqrt{P_1}\sin(\theta_1(t))$
($\omega_1(t) = \dot\theta_1(t)$ is the carrier frequency).
A local oscillator (LO - voltage-controlled oscillator VCO) signal is sinusoidal signal $\sqrt{P_2}\sin(\theta_2(t))$
with the frequency $\omega_2(t) = \dot\theta_2(t)$.
The block $90^o$ Hybrid combines inputs shifting phases by $90^o$ as follows:
\begin{equation}
\begin{aligned}
& E_1 =
		\frac{1}{2}m(t)\sqrt{P_1}\cos(\theta_1(t)) + \frac{1}{2}\sqrt{P_2}\cos(\theta_2(t)),
\\
& E_2 =
		\frac{1}{2}m(t)\sqrt{P_1}\cos(\theta_1(t)) - \frac{1}{2}\sqrt{P_2}\cos(\theta_2(t)),
\\
& E_3 =
		\frac{1}{2}m(t)\sqrt{P_1}\cos(\theta_1(t)) + \frac{1}{2}\sqrt{P_2}\cos(\theta_2(t) + \frac{\pi}{2}),
\\
& E_4 =
		\frac{1}{2}m(t)\sqrt{P_1}\cos(\theta_1(t)) - \frac{1}{2}\sqrt{P_2}\cos(\theta_2(t) + \frac{\pi}{2}).
\\
\end{aligned}
\end{equation}
Four outputs of the receivers are the following
\begin{equation}
\begin{aligned}
&
	I_1(t) =
	\frac{R}{8}
		\big(
			P_1 + P_2 + 2m(t)\sqrt{P_1P_2}\cos(\theta_1(t) - \theta_2(t))
	  	\big),
\\
&
	I_2(t)
	= \frac{R}{8}
		\big(
			P_1 + P_2 - 2m(t)\sqrt{P_1P_2}\cos(\theta_1(t) - \theta_2(t))
	  	\big),
\\
&
	I_3(t)
	= \frac{R}{8}
		\big(
		P_1 + P_2 + 2m(t)\sqrt{P_1P_2}\cos(\theta_1(t) - \theta_2(t) - {\pi\over2})
	  	\big),
\\
&
	I_4(t)
	= \frac{R}{8}
		\big(
	P_1 + P_2 - 2m(t)\sqrt{P_1P_2}\cos(\theta_1(t) - \theta_2(t) - {\pi\over2})
	  	\big).
\\
\end{aligned}
\end{equation}
The block TIA multiplies its inputs by certain factor $A$ and then
 subtracts two pairs of signals:
\begin{equation}
\begin{aligned}
&
	I_I(t) = AI_1(t) - AI_2(t)
  \\
  &
  \quad\quad=
	 {m(t)RA\sqrt{P_1P_2} \over 2}\cos(\theta_1(t) - \theta_2(t)),
\\
&
	I_Q(t) = AI_3(t) - AI_4(t)
  \\
  &
  \quad\quad\ =
	 {m(t)RA\sqrt{P_1P_2} \over 2}\cos(\theta_1(t) - \theta_2(t) - {\pi\over2}).
\\
\end{aligned}
\end{equation}
After multiplication the loop filter input becomes
\begin{equation}
\label{filter input}
\begin{aligned}
&
	\varphi(t) = I_I(t)I_Q(t) =
\\
&
	=
	{R^2A^2 P_1 P_2 \over 4}\cos(\theta_1(t)
	- \theta_2(t))\sin(\theta_1(t) - \theta_2(t))
\\
&
	=
	{R^2A^2 P_1 P_2 \over 8}\sin(2\theta_1(t) - 2\theta_2(t)).
\\
\end{aligned}
\end{equation}
It should be noted that unlike the classical PLL circuit
the filter input  of the considered model  (see \eqref{filter input})
depends on the difference of the phases only
(see also two-phase phase-locked loop model in \citep{BestKLYY-2014-IJAC,BestKLYY-2016}).
This allows one to derive a mathematical model of the optical Costas loop
in a signal's phase space
without additional engineering assumption
on the complete filtration of high-frequency signal component
\citep{LeonovKYY-2015-TCAS,BestKKLYY-2015-ACC}.

The relation between the input $\varphi(t)$
and the output $g(t)$ of the filter is as follows
\begin{equation}\label{loop-filter}
 \begin{aligned}
 & \dot x = A x + b \varphi(t),
 \ g(t) = c^*x + h\varphi(t),
 \end{aligned}
\end{equation}
where $A$ is a constant $n \times n$ matrix,
the vector $x(t) \in \mathbb{R}^n$ is a filter state,
$b,c$ are constant vectors,
and $x(0)$ is initial state of the filter.
The control signal $g(t)$ is used to adjust the VCO frequency to
the frequency of input carrier signal:
\begin{equation} \label{vco first}
   \dot\theta_2(t) = \omega_2(t) = \omega_2^{\text{free}} + K_{\rm VCO}g(t).
\end{equation}
Here $\omega_2^{\rm free}$ is a free-running frequency of the VCO
and $K_{\rm VCO}$ is the VCO gain.
Similarly, various nonlinear models of the VCO can be considered
(see, e.g., \citep{BianchiKLYY-2016}).
Note that the initial VCO frequency (at $t=0$) has the form
\begin{equation}
  \begin{aligned}
    &
    \omega_2(0) = \omega_2^{\rm free} + K_{\rm VCO}\alpha_0(0)
    +K_{\rm VCO}h\varphi(\theta_1(0) - \theta_2(0)).
  \end{aligned}
\end{equation}

\begin{figure}[H]
  \centering
  \includegraphics[width=\linewidth]{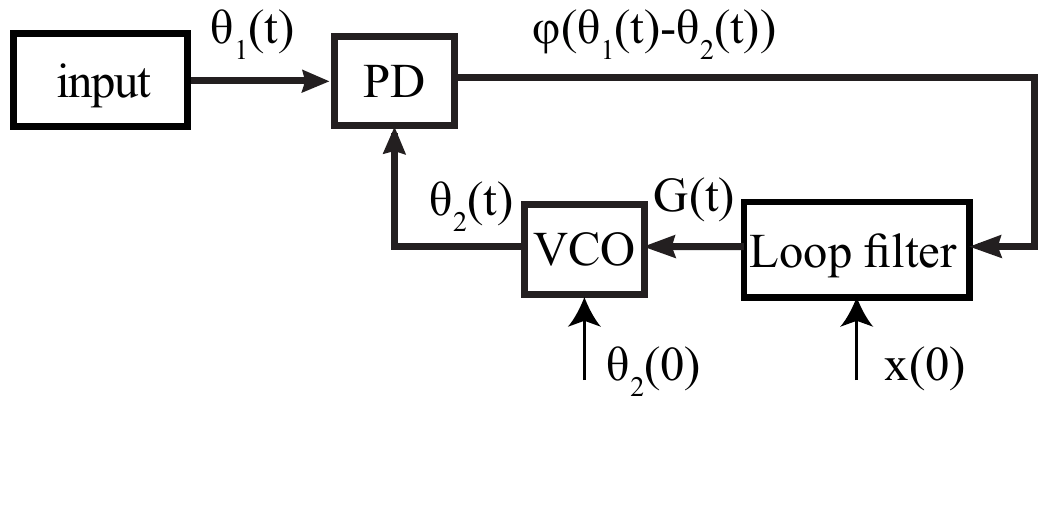}
  \caption{Optical BPSK Costas loop in signal's phase space.
  $\theta_2(0)$ is initial phase of the VCO, $x(0)$ is initial state of the filter.}
\label{phase model}
\end{figure}

Let  $\theta_{\Delta} = \theta_1 - \theta_2$.
If the frequency of the input carrier is constant:
\begin{equation}\label{omega1-const}
   \dot\theta_1(t) = \omega_1(t) \equiv \omega_1,
\end{equation}
then equations \eqref{loop-filter}-\eqref{vco first}
give the following mathematical model of the optical Costas loop in the signal's phase space:
\begin{equation} \label{mathmodel-class}
 \begin{aligned}
   & \dot x = A x + b \varphi(\theta_{\Delta}),\\
   & \dot \theta_{\Delta} =
   \omega_\Delta^{\rm free}-K_{\rm VCO}c^*x- h\varphi(\theta_{\Delta}),
 \end{aligned}
\end{equation}
where $|\omega_\Delta^{\rm free}| = |\omega_1 - \omega_2^{\rm free}|$ is called
a frequency deviation.
This model corresponds to the classical signal's phase model of the classical PLL
\citep{Viterbi-1966,Gardner-1966,ShahgildyanL-1966}, shown in Fig.~\ref{phase model}.

\section{The pull-in range computation}

One of the main engineering characteristics of the PLL-based circuits is the \emph{pull-in range} which corresponds to the global stability of the model.
The pull-in range is a widely used engineering concept
(see, e.g. \citep[p.40]{Gardner-1966}, \citep[p.61]{Best-2007}).
The following rigorous definition is from
\citep{KuznetsovLYY-2015-IFAC-Ranges,LeonovKYY-2015-TCAS,BestKLYY-2016}:
\emph{the largest interval of frequency deviations
$0 \leq |\omega_\Delta^{\rm{free}}|<\omega_{\text{pull-in}}$
such that the model in the signal's phase space
acquires lock for arbitrary initial phase difference and filter state
(i.e. any trajectory tends to a stationary point)
is called a \emph{pull-in range}, $\omega_{\text{pull-in}}$ is called
a \emph{pull-in frequency}.}

\subsection{Lead-lag filter} 
Consider the mathematical model of the optical Costas loop with lead-lag filter
in the signal's phase space.

\begin{figure}[ht]
  \centering\includegraphics[scale=0.4]{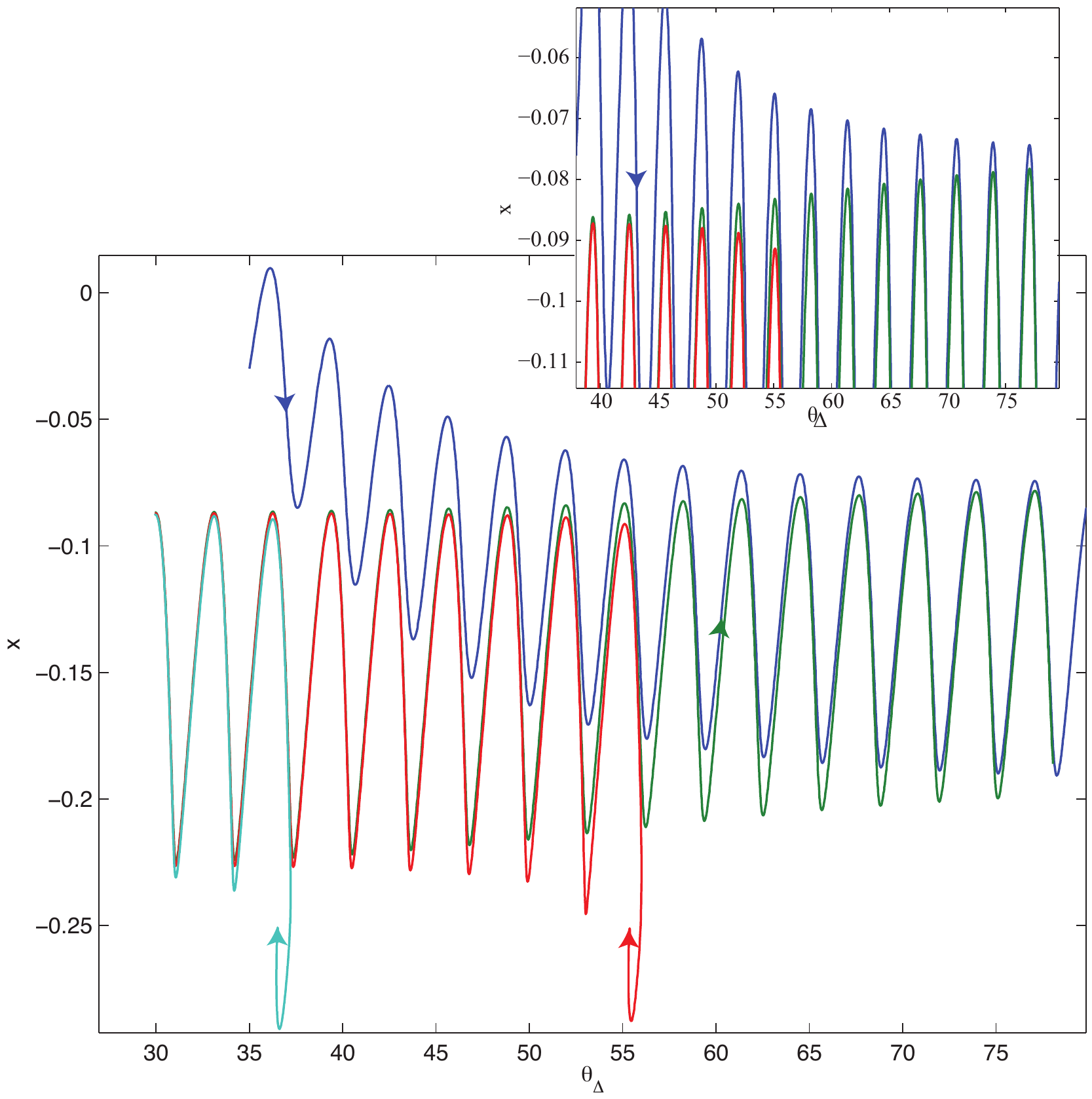}
  \caption{Phase portrait with stable and unstable cycles.}
  \label{semistable_phase}
\end{figure}
Fig.~\ref{semistable_phase} shows
phase portrait corresponding to
the signal's phase model  where
the carrier frequency is $10000$ rad/s,
the VCO free-running frequency is $10089.5$rad/s,
the lead-lag filter transfer function\footnote{
In engineering literature
filter transfer function is usually defined as $H(s) = c^*(sI-A)^{-1}b+h$ \citep{Best-2007},
at the same time in the control theory \citep{LeonovK-2014-book}
it is defined as $c^*(A-sI)^{-1}b-h$} is $H_{lf}(s) = \frac{a s+ \beta}{s + \alpha}$,
$\alpha = 63.1656$, $\beta = 63.1656$, $a = 0.2922$.
The solid blue line in Fig.~\ref{semistable_phase} corresponds
to the trajectory with the filter initial state
$x(0) = 0.6304$ and the VCO phase shift $0.3975$ rad.
This line tends to a periodic trajectory.
Thus in this case the model can not acquire lock.

The solid red line corresponds to the trajectory
with the filter initial state $x(0) = -0.1373$ and the VCO initial phase $24.3161$.
This trajectory lies just under an unstable periodic trajectory
(unstable limit cycle)
and tends to a stable equilibrium.
In this case the model acquires lock.

All the trajectories between the stable and unstable periodic trajectories
(see, e.g., a solid green curve)
tend to the stable one which is a stable limit cycle.
Therefore, if the gap between the stable and unstable periodic trajectories
is smaller than the discretization step,
then the numerical procedure may slip through the stable trajectory.
Thus, in this case we have coexistence of stable periodic trajectory
(which is a local hidden attractor\footnote{
Attractor is called a \emph{self-excited attractor} if its basin of attraction
intersects any arbitrarily small open neighborhood of an equilibrium,
otherwise it is called a \emph{hidden attractor}
\citep{LeonovKV-2011-PLA,LeonovKV-2012-PhysD,LeonovK-2013-IJBC,LeonovKM-2015-EPJST,Kuznetsov-2016}.
For self-excited attractors there is a transient process
from a small vicinity of an unstable equilibrium to an attractor,
which allows to find the attractor easily.
If  there is no such a transient process for an attractor,
it is called a hidden attractor.
For example, hidden attractors are attractors in systems
without equilibria or with only one stable equilibrium
(a special case of multistability and coexistence of attractors).
Some examples of hidden attractors can be found in
\cite{ShahzadPAJH-2015-HA,BrezetskyiDK-2015-HA,JafariSN-2015-HA,ZhusubaliyevMCM-2015-HA,SahaSRC-2015-HA,Semenov20151553,FengW-2015-HA,Li20151493,FengPW-2015-HA,Sprott20151409,PhamVVJ-2015-HA,VaidyanathanPV-2015-HA,ChenLYBXW-2015-HA,MenacerLC-2016-HA,Danca-2016-HA,Zelinka-2016-HA,DudkowskiJKKLP-2016,DancaKC-2016,KiselevaKL-2016-IFAC}.
}
--- locally attracting, closed, and bounded set in the cylindrical phase space)
and unstable periodic trajectory (repeller) close to each other;
the merge of these periodic trajectories leads to the birth of semistable trajectory
(semistable limit cycle) \citep{Gubar-1961,Shakhtarin-1969,BelyustinaBKS-1970,LeonovK-2013-IJBC,KuznetsovLYY-2014-IFAC}.
In this case numerical methods are limited
by the errors on account of the linear multistep integration methods
(see \citep{biggio2013reliable,biggio2014accurate}).
The corresponding difficulties of simulation in SPICE and MATLAB
are discussed in \citep{LeonovKYY-2015-TCAS,BianchiKLYY-2015,KuznetsovLYY-2017-CNSNS}.

\begin{figure}[H]
	\includegraphics[width=\linewidth]{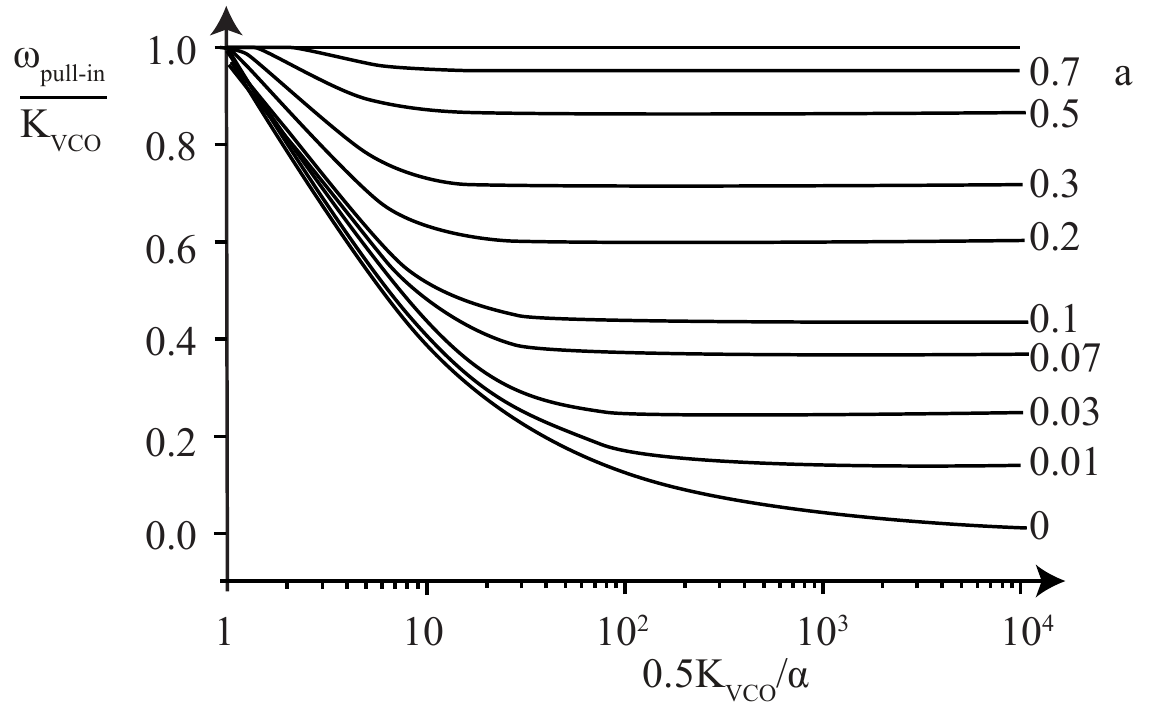}
	\caption{Pull-in range diagram for the optical Costas loop}
	\label{optical-pull-in}
\end{figure}

The above example shows that, in general, the pull-in of the model
with a lead-lag filter is bounded
and demonstrates the difficulties of numerical estimation of the pull-in range.
In this case the pull-in range estimation can be done
by analytical-numerical methods based on phase-plane analysis,
however, it is a challenging task due to hidden oscillations
(see, e.g. \citep{Shakhtarin-1969,BelyustinaBKS-1970,Shalfeev-2013-book,LeonovK-2013-IJBC}).
Corresponding diagram for pull-in range is on Fig.~\ref{optical-pull-in}.

\begin{figure}[ht]
  \centering
  \includegraphics[width=0.95\linewidth]{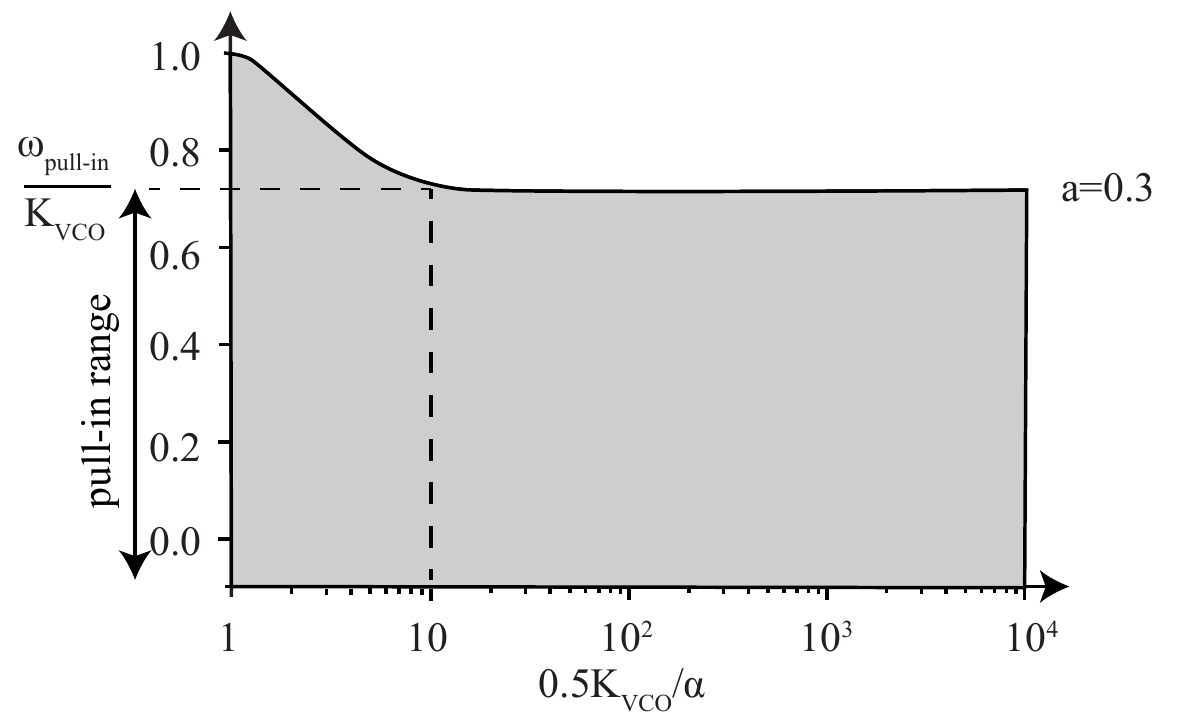}
  \caption{Pull-in range diagram explanation.}
  \label{bif_diagram_explain}
\end{figure}

Consider, e.g., parameters $a = 0.3$ and
 corresponding curve in Fig.~\ref{bif_diagram_explain}.
Here normalized pull-in frequency $\frac{\omega_\Delta^{\rm{free}}}{K_{\rm VCO}}$
is plotted for various values of the loop gain $K_{\rm VCO}$.
For all points on the diagram in Fig.~\ref{bif_diagram_explain} below the curve (filled with grey) system \eqref{mathmodel-class}
is globally asymptotically stable.
Pull-in frequency diagram for different parameters of the loop filter is
in Fig.~\ref{optical-pull-in} for sinusoidal PD characteristic.

\subsection{PI filter}
The filter -- perfect integrator is preferred
(if it can be implemented in the considered architecture)
since it allows us to prove analytically that the pull-in range is infinite.

Since PLL-based circuits are nonlinear control systems
for their global analysis it is essential to apply
the stability criteria, which are developed in control theory,
however their direct application to analysis of the PLL-based models
is often impossible, because such criteria are usually not adapted for
the cylindrical phase space\footnote{For example,
in the classical Krasovskii--LaSalle principle on global stability
the Lyapunov function has to be radially unbounded
(e.g. $V(x,\theta_\Delta) \to +\infty$ as $||(x,\theta_\Delta)|| \to +\infty$).
While for the application of this principle to the analysis of phase synchronization systems
there are usually used Lyapunov functions periodic in $\theta_\Delta$
(e.g. $V(x,\theta_\Delta)$  is bounded
for any $||(0,\theta_\Delta)|| \to +\infty$),
and the discussion of this gap is often omitted
(see, e.g. patent \citep{Abramovitch-2004} and works \citep{Bakaev-1963,Abramovitch-1990,Abramovitch-2003}).
Rigorous discussion can be found, e.g. in \citep{GeligLY-1978,LeonovK-2014-book}.
};
The corresponding modifications of classical stability criteria
for the nonlinear analysis of control systems in cylindrical phase space
were developed in the second half of the 20th century
(see, e.g. \citep{GeligLY-1978,LeonovRS-1992,LeonovKS-2009,LeonovK-2014-book}).

\begin{figure}[h]
\centering
  \includegraphics[width=\linewidth]{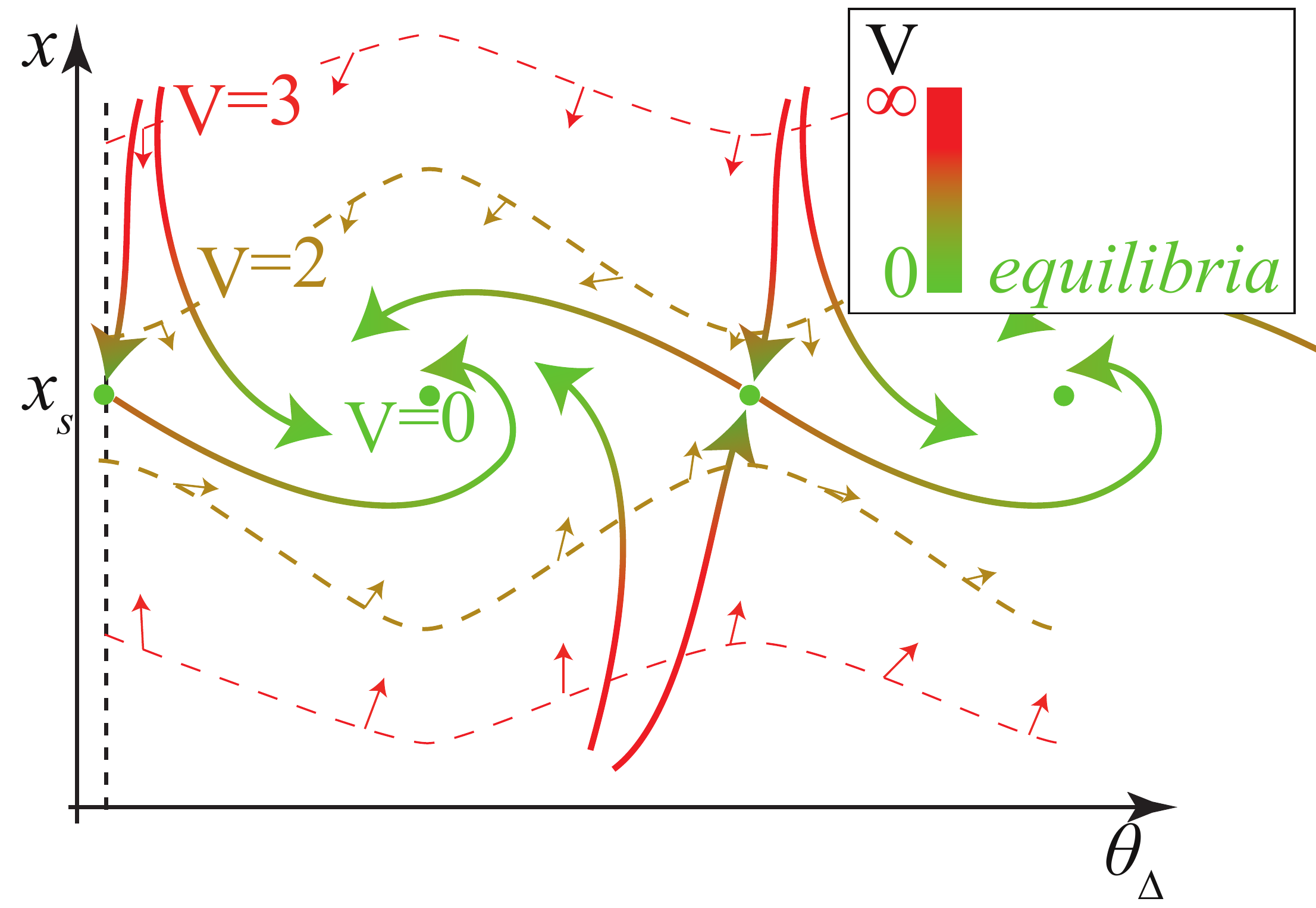}
  \caption{Lyapunov function is nonincreasing along the trajectories,
   periodic case.}
  \label{lyapunov-periodic-expl}
\end{figure}



Consider model \eqref{mathmodel-class}
with a proportionally-integrating (active PI) filter
\begin{equation}
\label{perfect-pi-pll}
  \begin{aligned}
    & \dot x = \frac{1}{\tau_1}\sin(2\theta_\Delta),\\
    & \dot\theta_\Delta =
    \omega_\Delta^{\rm{free}}
    - K_{\rm VCO}(x + \frac{\tau_2}{\tau_1}\sin(2\theta_\Delta)),
  \end{aligned}
\end{equation}
where
\begin{equation}
	\begin{aligned}
		& \tau_1 > 0,\ \tau_2 > 0,
	\end{aligned}
\end{equation}
and the following Lyapunov function \citep{Bakaev-1965,LeonovKYY-2015-TCAS}
\begin{equation}
  \begin{aligned}
  &V(x, \theta_\Delta) =
  \frac{K_{\rm VCO}\tau_1}{2} \left(x - \frac{\omega_\Delta^{\rm free}}{K_{\rm VCO}}\right)^2
  + \int \limits_{0}^{\theta_\Delta} \sin(2s) ds.\\
  \end{aligned}
\end{equation}
For its derivative we get (see Fig.~\ref{lyapunov-periodic-expl})
\begin{equation}
\label{pi-pull-in-cond}
\begin{aligned}
  &\dot{V}(x, \theta_\Delta) =
  \frac{1}{\tau_1}\sin(2\theta_\Delta)
  \frac{K_{\rm VCO}\tau_1}{2}
  2
    \left(
      x - \frac{\omega_\Delta^{\rm free}}{K_{\rm VCO}}
    \right)\\
  &+\sin(2\theta_\Delta)
  \left(
    \omega_\Delta^{\rm{free}} - K_{\rm VCO}(x + \frac{\tau_2}{\tau_1}\sin(2\theta_\Delta))
  \right)\\
  &= -\frac{K_{\rm VCO} \tau_2}{\tau_1}^2\sin^2(2\theta_\Delta) \leq 0.
\end{aligned}
\end{equation}
Remark that for any $\omega_\Delta^{\text{free}}$
the set $\dot V(x,\theta_\Delta)\equiv0$ does not contain the whole trajectories except for equilibria. 
Thus, form the modification of classical direct Lyapunov method
for cylindrical phase space
(\citep{GeligLY-1978,LeonovK-2014-book}
it follows that the pull-in range is infinite
since \eqref{pi-pull-in-cond} holds for arbitrary $\omega_\Delta^{\rm free}$.

\section{Conclusion}
In this work we consider
a mathematical model of the optical Costas loop
and study its pull-in and lock-in ranges.
For the lead-lag filter we show that numerical estimation of the
pull-in range is a challenging task due to the existence of hidden oscillations.
For the active PI filters we prove that the pull-in range is infinite
by application of special modification of the direct Lyapunov method for the cylindrical phase space.

\section*{Acknowledgment}
This work was supported by the Russian Science Foundation (project 14-21-00041).
The authors would like to thank Roland~E.~Best
(the author of the bestseller on PLL-based circuits \citep{Best-2007})
for valuable discussion on the stability ranges,
and Professor Prof. Werner Rosenkranz (Kiel University)
for demonstration of the optic Costas loop operation.

\bibliographystyle{ifacconf}

\end{document}